\documentclass[12pt]{article} 
\usepackage{amsmath, amsthm}
\input amssym.def
\input amssym
\def \End{{\rm End}}

\def \<{\langle} 
\def \>{\rangle}

\def \pf {\noindent {\bf Proof:} \,}

\newcommand{\NO}{\,{\raise0.25em\hbox{$\mathop{\hphantom{\cdot}}%
\limits^{_{\circ}}_{^{\circ}}$}}\,}

\theoremstyle{plain}

\begin{document}

\newtheorem{thmm}{Theorem}
\newtheorem{remm}[thmm]{Remark}
\newtheorem{thm}{Theorem}[section]
\newtheorem{prop}[thm]{Proposition}
\newtheorem{coro}[thm]{Corollary}
\newtheorem{lem}[thm]{Lemma}
\newtheorem{rem}[thm]{Remark}
\newtheorem{de}[thmm]{Definition}


\begin{center}
{\Large {\bf Local and Semilocal Vertex Operator Algebras  }} \\
\vspace{0.5cm} Chongying Dong\footnote{Partially supported by NSF
grants and Research grants from UC Santa Cruz} and Geoffrey
Mason\footnote{Partially supported by NSF}\\ Department of Mathematics
University of California, Santa Cruz, 95064
\end{center}

\begin{abstract} We initiate a general structure theory for
vertex operator algebras $V$. We discuss the
center and the blocks of $V$,
  the Jacobson radical and solvable radical, and
   local vertex operator algebras. The main consequence
of our structure theory is that if $V$ satisfies some mild conditions,
  then it is necessarily semilocal, i.e. a direct
sum of local vertex operator algebras.
\end{abstract}

\section{ Introduction}
\setcounter{equation}{0}

The purpose of this paper is to establish some results concerning the
algebraic structure of vertex operator algebras. Much
of the existing literature is devoted to the study of
\emph{simple} vertex operator algebras, which are certainly very important.
However, we perceive the need to have available
  a general structure theory
for vertex operator algebras. The results in the present paper
may serve as a step towards this goal.

It is a well-known heuristic,  a consequence of the
commutivity and associativity axioms [LL], that vertex operator
algebras behave, in some ways,
like finite-dimensional commutative associative algebras.
  It is therefore natural to see to what extent one
can develop a structure theory for vertex operator algebras
that parallels that of finite-dimensional algebras,
  which is particularly transparent if the algebra
is commutative and associative.
This entails the idea of the \emph{center}
and various types of \emph{radical ideals} of $V$.
  The center of $V$ is defined as follows:
\begin{equation}\label{eq: Z(V)}
Z(V) = \{ v \in V \ | \ L(-1)v = 0 \}.
\end{equation}
It is both a finite-dimensional commutative, associative algebra
with respect to the $-1$ product in $V$,
and a vertex subalgebra of $V$. The center provides the connection
between vertex operator algebras
and the classical theory of algebras. We define $J(V)$ along the lines of
one of the definitions for associative algebras (some of the other definitions
do not work so well for vertex operator algebras). Thus,
\begin{equation}\label{eq: J(V)}
J(V) = \cap \ M
\end{equation}
where the intersection ranges over the maximal ideals of $V$.
For obvious reasons we call it the \emph{Jacobson radical} of $V$.
It is not hard to see that $J(V)$ is the smallest ideal of $V$
such that the quotient vertex operator
algebra $V/J(V)$ is a semisimple vertex operator algebra. Furthermore 
it is always true
(Theorem 3.6) that
$J(V) \cap Z(V) = J(Z(V))$.
We also introduce nilpotent and solvable ideals
and the \emph{solvable radical} of a vertex operator algebra.

We consider the decomposition of
  $V$ into \emph{blocks}, i.e. the
indecomposable direct summands of the adjoint module. These are
vertex operator algebras themselves, having the same
central charge as $V$. An important point
(Proposition 2.6) is that a vertex operator algebra
is indecomposable precisely when its center is a local algebra.
These general results lead to the main structural results
of the paper. To describe them, let us define a vertex operator algebra to be
\emph{local} in case it has a unique maximal ideal.
In this case, $J(V)$ is the maximal ideal and
$V/J(V)$ is a simple vertex operator algebra. We say that $V$ is 
\emph{semilocal}
if $V$ is (isomorphic to) a direct sum of local vertex operator 
algebras (which necessarily
have the same central charge).

\begin{thmm}\label{t1}  A vertex operator algebra $V$ is semilocal
under either of the following two conditions: \\
(i) There are no nonzero weight spaces $V_n$ for $n<-1$; \\
(ii) The Jacobson radical of $V$ coincides with the solvable
radical.
\end{thmm}

Our approach to the question of semilocality and the structure of
local vertex operator algebras under the hypotheses of part (i)
of Theorem 1
  depend on the nature of the nonassociative algebra structure on the
zero weight space $V_0$ induced by the $-1$st and related
products. We shall establish the following result, which readily implies
  part (i) of Theorem 1.

\begin{thmm}\label{t2} Suppose that $V$ is a vertex operator algebra such that
$V_n = 0$ for $n < -1$. Then the following are equivalent:\\
(a) $V$ is  local;\\
(b) $V$ is indecomposable;\\
(c) $Z(V)$ is a local algebra; \\
(d) $V_0$ is a commutative, power associative local algebra with respect to the
product $a*b = 1/2(a(-1)b + b(-1)a)$.
\end{thmm}

\begin{remm}{\rm  (i) We call a power associative algebra $A$  \emph{local}
  in case the
nil radical (largest ideal for which all elements are nilpotent) has 
codimension $1$
  and is the unique maximal ideal in $A$. (See Section 4 for more details.) \\
(ii) If $V$ satisfies the stronger condition that $V_n = 0$ for $n<0$ then
$V_0$ is actually a commutative, associative local algebra and the $*$ product
in (d) is the usual product $a*b = a(-1)b$. In this case, (d) says that
$V_0$ is a local algebra in the usual sense.}
\end{remm}

Our approach to part (ii) of Theorem 1 might appear somewhat different
to that of part (i), depending as it does
on general properties of the radical ideals. But the goal in each case is to
establish a sort of idempotent lifting theorem (cf. Section 4), and the two
approaches eventually converge. This raises the general question:
  \begin{eqnarray*} \mbox{are \emph{all} vertex operator algebras semilocal? }
\end{eqnarray*}
This amounts
to asking if conditions (a) and (b) in Theorem 2 are
  always equivalent.

The paper is organized as follows: in Section 2 we develop the ideas
of the center and the blocks of a vertex operator algebra $V$. Our
development here is related to parts of Li's paper [L]. Indeed the center
of $V$, being the kernel of $L(-1)$, is in a sense dual to the cokernel
of $L(1)$, which is the focus of Li's paper. In Section 3 we consider 
various radical ideals,
in particular the Jacobson radical of $V$, while in
Section 4 we prove the main Theorems. In the final Section 5
we discuss some examples.

\section{The center and the blocks of a vertex operator algebra}
\setcounter{equation}{0}

We begin with some elementary results about vertex operator algebras $V$.
  We consider $V$ as a $Z$-graded space
with respect to the usual $L(0)$-grading,
\begin{equation}\label{eq: vgrading}
V = \bigoplus_{n \in Z} V_n.
\end{equation}
  As in (\ref {eq: Z(V)}), $Z(V)$ denotes the center of $V$.

\begin{lem}\label{l2.1}
Let $v \in V$. Then the following are equivalent: \\
(a) $v \in Z(V)$; \\
(b) The vertex operator for $v$ is a constant, i.e. $Y(v, z) = v(-1)$.
\end{lem}
\pf  If (a) holds then $L(-1)v = 0$, and therefore
$(L(-1)v)(n) = 0$ for all $n$. Since
\begin{eqnarray*}
(L(-1)v)(n) = -nv(n-1)
\end{eqnarray*}
it follows that $v(n-1) = 0$ for $n \neq 0$. Hence, $Y(v, z) = \sum 
v(n)z^{-n-1} = v(-1)$.
This shows that (a) $\Rightarrow$  (b). If (b) holds then we have
\begin{eqnarray*}
0 = \frac{d}{dz} Y(v, z) = Y(L(-1)v, z),
\end{eqnarray*}
so that
\begin{eqnarray*}
0 = Y(L(-1)v, z)\mathbf{1} = L(-1)v + O(z).
\end{eqnarray*}
Therefore (a) holds. \qed

\begin{lem}\label{l2.2} The following hold: \\
(a) $Z(V)$ is annihilated by all Virasoro operators $L(n)$ with $n \geq -1$. \\
(b) $Z(V) \subseteq V_0$ consists of primary states.
\end{lem}

\pf (a) Let $v \in Z(V)$. For some positive integer $k$ we
have $L(k)v = 0$. Then
\begin{eqnarray*}
0 = [L(k), L(-1)]v = (k+1)L(k-1)v,
\end{eqnarray*}
so that also $L(k-1)v = 0$. From this we deduce that $L(n)v = 0$ for 
all $n\geq 0$, and
in particular (a) holds. Part (b) is a special case of (a),
being equivalent to the equalities $L(2)v = L(1)v = L(0)v = 0$ for $v 
\in Z(V)$. \qed

\begin{lem}\label{l2.3} The $(-1)$st product $ab = a(-1)b$ for $a, b \in 
Z(V)$ gives $Z(V)$
the structure of a commutative, associative algebra with identity $\mathbf{1}$.
\end{lem}

\pf  Let $a, b \in Z(V)$. Then
\begin{eqnarray*}
L(-1)a(-1)b = [L(-1), a(-1)]b = (L(-1)a)(-1)b = 0.
\end{eqnarray*}
This shows that $a(-1)b  \in Z(V)$,
  so that $Z(V)$ is closed with respect to the product $a(-1)b$.
 From skew-commutativity we get
\begin{eqnarray*}
b(-1)a = \sum _{n \geq 0} (-1)^nL(-1)^na(n-1)b = a(-1)b,
\end{eqnarray*}
where we have used Lemma 2.1. So the product is commutative. A
similar proof (which we omit) using associativity shows that the 
product is also
associative, and the Lemma follows. \qed

For an element $v \in Z(V)$ we define the \emph{annihilator} of $v$ as
\begin{eqnarray*}
\mbox{Ann}_V(v) = \mbox{Ann}(v) = \{ u \in V | v(-1)u = 0 \}.
\end{eqnarray*}

\begin{lem}\label{l2.4}
The following hold: \\
(a) $[v(-1), Y(u, z)] = 0$ for any $v \in Z(V), u \in V$ ; \\
(b) Ann$(v)$ is an ideal in $V$ for any $v \in Z(V)$.
\end{lem}

\pf Let $v \in Z(V)$ and $u \in V$. By the commutativity 
formula we have
\begin{eqnarray*}
[v(-1), u(n)] = \sum_{m\geq 0}(-1)^m (v(m)u)(m+n-1) = 0,
\end{eqnarray*}
   the second equality following from Lemma 2.1. So (a) holds.
  As for (b), suppose that $a \in$ Ann$(v)$. Using part (a) we get
\begin{eqnarray*}
v(-1)u(n)a = u(n)v(-1)a = 0.
\end{eqnarray*}
This shows that Ann$_V(v)$ is closed under the (right) action of all 
modes $u(n)$,
and hence is an ideal of $V$.  \qed

We define the \emph{endomorphism ring} End($V$) = Hom$_V(V, V)$ in
the usual way, namely
\begin{eqnarray*}
\mbox{End}(V) = \{ \varphi \in \mbox{End}_{\mathbf{C}}(V, V) \ |  \ 
\varphi Y(u, z)
= Y(u, z) \varphi, \mbox{all} \ u \in V \}.
\end{eqnarray*}
For $a \in Z(V)$ we define a linear map $\varphi_a: V \rightarrow V$ via
\begin{equation}\label{eq: phimap}
\varphi_a: v \mapsto a(-1)v.
\end{equation}
   Note that as a direct consequence of Lemma 3.4(a), we have
\begin{equation}\label{eq: varphicontain}
\varphi_a \in \mbox{End}(V).
\end{equation}

\begin{prop}\label{p2.5}
The map $a \mapsto \varphi_a$ induces a 
natural isomorphism of rings
\begin{eqnarray*}
Z(V)   \stackrel{\cong}{\longrightarrow} \mbox{End}(V).
\end{eqnarray*}
\end{prop}

\pf  Let $\varphi \in$ End$(V)$ and $u \in V$. Because $\varphi$
commutes with all modes $u(n)$ then
\begin{eqnarray*}
\varphi(u) = \varphi(u(-1)\mathbf{1}) = u(-1) \varphi(\mathbf{1}).
\end{eqnarray*}
This shows that $\varphi$ is uniquely determined by the image of the vacuum
state. Furthermore
\begin{eqnarray*}
L(-1) \varphi(\mathbf{1}) = \varphi L(-1) \mathbf{1} =  0,
\end{eqnarray*}
so that $\varphi (\mathbf{1}) \in Z(V)$. On the other hand, for $a \in Z(V)$
we have $\varphi_a(\mathbf{1}) = a$, as follows from (\ref {eq: phimap}). This
shows that the map $a \mapsto \varphi_a$ induces a linear isomorphism
from $Z(V)$ to End($V$). That it is an isomorphism of rings follows from the 
identity
$(a(-1)b)(-1) = a(-1)b(-1)$ for $a, b \in Z(V)$. \qed

We will refer to the idempotent elements of $Z(V)$
as the \emph{idempotents of $V$}.
We similarly transfer other standard constructions and language
  concerning idempotents
(for example \emph{orthogonal} idempotents and \emph{primitive} idempotents)
from $Z(V)$ to $V$. Thus $V$ has a \emph{unique} set of primitive idempotents
$e_1, ..., e_n$, and we have
\begin{equation}\label{eq: 1decomp}
e_1+ ... + e_n = \mathbf{1}.
\end{equation}
For $a \in Z(V)$, both image and kernel of $\varphi_a$ are ideals in $V$,
indeed the kernel is nothing but the annihilating ideal Ann$(a)$
(cf. Lemma 2.4(b)). Of course, in general we have Im$\varphi_a \cong 
V/$Ann$(a)$.
But if $e$ is an idempotent of $V$ then there is a splitting
\begin{equation}\label{eq: splitting}
V = e(-1)V \oplus (\mathbf{1} - e)(-1)V = \ \mbox{Im} \ \varphi_e 
\oplus \mbox{ker} \ \varphi_e.
\end{equation}
Now set $V^i = e_i(-1)V, 1 \leq i \leq n$. We therefore obtain a 
decomposition of
$V$ into ideals
\begin{equation}\label{eq: Vdecomp1}
V = V^1 \oplus ... \oplus V^n
\end{equation}
such that $e_i = \mathbf{1}_i$ is the vacuum state for $V^i$. The uniqueness
of the primitive idempotents for $V$ entails the uniqueness of this
decomposition. We call the ideals $V^i$ the \emph{blocks} of $V$, and
(\ref {eq: Vdecomp1}) the decomposition of $V$ into blocks.

\begin{prop}\label{p2.6} The following are equivalent: \\
(a) $V$ is indecomposable as a $V$-module. \\
(b) $\End(V) = Z(V)$ is a local commutative, associative algebra.
\end{prop}

\pf  We already know (Proposition 2.5) that End$(V) = Z(V)$
is always a commutative, associative algebra with identity element 
$\mathbf{1}$.
Then the equivalence of (a) and (b) has essentially the same proof
as for finite-dimensional (associative) algebras. Namely, every idempotent
in End$(V)$ determines a splitting (\ref {eq: splitting}) into ideals 
of $V$ as above.
Conversely if $V = A \oplus B$ is a splitting into ideals of $V$ then
projection onto $A$ is an idempotent in End$(V)$. We leave further details
to the reader. \qed.

We have proved most of the next result.

\begin{thm}\label{t2.7}
(Block decomposition of $V$)  A vertex operator algebra
$V$ has a unique decomposition (\ref {eq: Vdecomp1}) into blocks 
$V^1, ..., V^n$.
The blocks of $V$ are the indecomposable direct summands of the 
adjoint module $V$.
The number of blocks $n$ is equal to the number of primitive 
idempotents in $Z(V)$.
\end{thm}

\pf It only remains to explain why the blocks are indecomposable. Indeed, 
it is clear that
\begin{eqnarray*}
Z(V) = Z(V^1) \oplus ... \oplus Z(V^n),
\end{eqnarray*}
and that $e_i$ is the identity element of $Z(V^i)$. Thus, $e_i$ is the
\emph{unique} nonzero idempotent in $Z(V^i)$. So $Z(V^i)$ is a local algebra,
and therefore $V^i$ is indecomposable by Proposition 2.6. This completes the
proof of Theorem 2.7.

Indecomposable vertex operator algebras arise naturally in
representation theory. First some standard definitions: given a 
vertex operator algebra $V$
and a (nonzero) $V$-module $(M, Y_M)$, the annihilator of $M$ is
\begin{eqnarray*}
\mbox{Ann}(M) = \{ v \in V \ | \ Y_M(v, z)M = 0 \}.
\end{eqnarray*}
It is an ideal of $V. \ M$ is called \emph{faithful} if
Ann$M = 0$. So in general, $M$ is a faithful module
over $V/$Ann($M)$.

\begin{lem}\label{l2.8} Suppose that $M$ is a faithful, indecomposable $V$-module.
Then the following hold: \\
(a) $V$ is indecomposable;\\
(b) If $M$ is simple then $Z(V) = \mathbf{C1}$.
\end{lem}

\pf By a standard argument, since $M$ is indecomposable
then \linebreak
Hom$_V(M,M)$ is a local algebra. Because of Lemma 2.4(a)
there is a natural injection of rings $Z(V) \rightarrow$ Hom$_V(M, M)$,
and therefore $Z(V)$ is also a local algebra. Now part (a) follows
from Proposition 2.6. If $M$ is simple then
Hom$_V(M, M) = \mathbf{C}$ by Schur's Lemma, so (b) holds. \qed

\section{Radical ideals of a vertex operator algebra.}
\setcounter{equation}{0}

We discuss various types of radical ideals
for a vertex operator algebra. They are modelled
in a rather obvious way on corresponding ideas from
the theory of algebras, both associative and nonassociative.
We repeatedly use various elementary properties
of ideals $I \subseteq V$ (cf. [LL]): left ideals are necessarily $2$-sided,
and they are $Z$-graded subspaces of $V$. That is,
\begin{eqnarray*}
I = \bigoplus_{n \in Z} I_n, \     I_n = I \cap V_n.
\end{eqnarray*}

It is evident that there is a \emph{unique} ideal $T = T(V) \subseteq 
V$ maximal with
respect to having the property that the smallest weight of
a nonzero state in $T$ is at least $2$. (By Lemma 2.2(b), this is the same as
requiring that $T_0 = T_1 = 0$.)
Note also that $T(V/T(V)) = 0$, so that $T(V)$ is a
type of radical, though of a very simple kind.
We call $T$ the \emph{trivial radical} because
it has little effect on the main issues we are concerned with, which 
center around
the weight-space $V_0$. Here is a simple example of this idea:

\begin{lem}\label{l3.1}
Let $T$ be the trivial radical of the vertex 
operator algebra $V$.
  Then the following are equivalent:\\
(a) $V$ is indecomposable; \\
(b) $V/T$ is indecomposable.
\end{lem}

\pf Let $a \in V$ be such that $a+T \in Z(V/T)$. Without loss we
may choose $a \in V_0$. Then
$L(-1)(a+T) = T$, so that $L(-1)a \in T \cap V_1 = 0$. This shows that
we must have $a \in Z(V)$. Because $T_0 = 0$ it follows that the
projection $V \rightarrow V/T$ restricts to an isomorphism
$Z(V) \stackrel{\cong}{\rightarrow} Z(V/T)$. Now the Lemma follows 
immediately from
Proposition 2.6. \qed

\begin{lem}\label{l3.2} Let $V$ be a vertex operator algebra such that
the trivial radical $T(V)$ vanishes. Then every nonzero ideal of
$V$ contains a \emph{minimal} ideal (that is, a nonzero irreducible
$V$-submodule).
\end{lem}

\pf  Let $I$ be a nonzero ideal of $V$. As $T(V) = 0$ then
$I_0+I_1 \neq 0$. Since $I_0 + I_1$ has finite dimension it is clear
that we can find a nonzero ideal $P$ of $V$ contained in $I$
for which $P_0+P_1$ has minimal dimension. Moreover, we
may assume that $P$ is generated (as $V$-module) by $P_0+P_1$.
We assert that $P$ is a minimal ideal. Indeed, if $Q \subseteq P$
is a non-zero ideal of $V$ then $0 \neq Q_0+Q_1 \subseteq P_0+P_1$.
Choice of $P$ forces  $Q_0+Q_1 = P_0+P_1$, so that $P = Q$
since $P$ is generated by $P_0+P_1$. This completes the proof of
the Lemma. \qed

We next consider the Jacobson radical $J(V)$.
  Let $\mathcal{M}$ be the set of maximal (proper) ideals of $V$.
It is evident that $\mathcal{M}$ is non-empty, indeed that
  every (proper) ideal of $V$ is contained in a maximal ideal.
Moreover, an ideal $M$ is maximal if, and only if, $V/M$ is a simple 
vertex operator algebra.
We then define $J(V)$
   as in (\ref{eq: J(V)}).
\begin{lem}\label{l3.3} Suppose that $I$ is an ideal in $V$ with $I_0 = 
\{0\}$. Then $I \subseteq J(V)$.
\end{lem}

\pf If there is
a maximal ideal $M$ in $V$ which does not contain $I$ then
$V = M + I$. Then $\mathbf{1} \in V_0 = I_0 + M_0 = M_0 \subseteq M$, 
contradiction.
\qed

Lemma 3.3 implies that $T(V) \subseteq J(V)$.

\begin{prop}\label{p3.4} The following hold: \\
  (a) If $I $ is an ideal in $V$ such that $V/I$ is semisimple, then 
$J(V) \subseteq I$; \\
(b) $V/J(V)$ is semisimple.
\end{prop}

\pf Part (a) follows by a standard argument which we omit.
As for part (b), let $J = J(V)$ and let the maximal ideals of $V$ be
$M^1, M^2, ..., M^k, ...$. A standard argument shows that each intersection
$D^k = \bigcap_{i=1}^k M^i$
is such that $V/D^k$ is semisimple, so it
is enough to show that there is an integer $n$ such that  $D^n = J$.
Because $V_0$ has finite dimension there is certainly an integer $n$
such that $(D^n)_0 = J_0$.   Then in the quotient vertex operator algebra
$V/J$ we have $(D^n/J)_0 = 0$. By Lemma 3.3,
\begin{eqnarray*}
  D^n/J \subseteq J(V/J) = \{0\}.
\end{eqnarray*}
Thus using part(a) we get $J \subseteq D^n \subseteq J$. The Lemma is 
proved.  \qed

It follows immediately from Proposition 3.4 that we have

\begin{coro}\label{l3.5} Let $V$ be a vertex operator algebra. Then $V$ 
has only a finite number
of distinct maximal ideals, say $M^1, ..., M^k$. Moreover $V/J(V)$ is 
semisimple and
  has a decomposition
\begin{equation}\label{moredecomp}
V/J(V) = (V^1/J(V)) \oplus ... \oplus (V^k/J(V))
\end{equation}
into the direct sum of $k$ simple vertex operator algebras $V^i/J(V)$. 
\end{coro}
\begin{lem}\label{l3.6}
Suppose that $V$ is a
  direct sum of $k$ simple vertex operator algebras
$V^1, ... , V^k$. Let $\bf{1}$$_i$ be the vacuum vector of $V^i$. Then
\begin{eqnarray*}
Z(V) = \bf{C}\bf{1}_{\mbox{1}} \oplus ... \oplus \bf{C}\bf{1}_{\mbox{k}}
\end{eqnarray*}
is a semisimple
algebra of rank $k$ and $\{\bf{1}$$_1, ...,  \bf{1}$$_k\}$
is the complete set of primitive idempotents of $Z(V)$.
\end{lem}

\pf  It is clear that $Z(V) = Z(V^1) \oplus ... \oplus 
Z(V^n)$. Now the present Lemma follows
from Lemma 2.8(b). \qed

\begin{thm}\label{t3.7} For any vertex operator algebra $V$, we have
\begin{eqnarray*}
J(V) \cap Z(V) = J(Z(V)).
\end{eqnarray*}
\end{thm}
\pf  By Theorem 2.7 we may assume that $V$
is indecomposable, and therefore ( Proposition 2.6) that
$Z = Z(V)$ is a local commutative algebra.
Since $V/J$ is semisimple
(Proposition 3.4(b)) then
$Z(V/J)$ is semisimple by Lemma 3.6. Being a subalgebra of $Z(V/J)$, 
it follows that
$Z+J/J \cong Z/Z \cap J$ is also semisimple, and hence that $J(Z) 
\subseteq Z \cap J$.
  Since $J(Z)$ has codimension $1$ in $Z$, and since $\mathbf{1} \notin J$,
the Theorem follows. \qed

It is possible to refine Theorem 3.7. In order to do this, we make 
the following definitions.
Let $t \geq 1$ be an integer, and $I$ an ideal in a vertex operator 
algebra $V$. Set
\begin{eqnarray*}
I^t = \mbox{subspace of $V$ spanned by all states of the form} \\
   a_1(n_1)...a_{t-1}(n_{t-1})a_t, \ \mbox{all} \  a_1, ..., a_t \in 
I, n_1, ..., n_t \in Z.
\end{eqnarray*}
  It is easy to see that $I^t$ is an ideal of $V$ whenever $I$ is. In this
way we obtain two descending sequences of ideals
\begin{eqnarray*}
I \supseteq I^2 \supseteq I^3 \supseteq ...\supseteq I^r \supseteq ...
\end{eqnarray*}
and
\begin{eqnarray*}
I \supseteq I^2 \supseteq (I^2)^2 \supseteq ...\supseteq (I^{(r)})^2 
\supseteq ...
\end{eqnarray*}
where we inductively define $I^1 = I$ and $I^{(r)} = (I^{(r-1)})^2$.
We say that $I$ is \emph{nilpotent} in case there is $r \geq 1$ such 
that $I^r = 0$,
and \emph{solvable} if there is $r \geq 1$ such that $I^{(r)} = 0$. 
Observe that
$I^{(r)} \subseteq I^{2^r}$, so that $I$ nilpotent $\Rightarrow I$ solvable.

\begin{thm}\label{t3.8} The following hold for a vertex operator algebra $V$: \\
(a) $V$ has a \emph{unique} maximal nilpotent ideal $N(V)$; \\
(b) $N(V) \subseteq J(V)$; \\
(c) $Z(V) \cap N(V) = J(Z(V))$.
\end{thm}
\pf  To prove (a), it is enough to show that the sum $I^1 + I^2$
of two nilpotent ideals $I^1, I^2 \subseteq V$ is again nilpotent. Let $t$
be a positive integer such that $(I^j)^t = 0, j = 1,2$, and choose elements
$x_1, ..., x_{2t} \in I^1, y_1, ..., y_{2t}\in I^2$. It suffices to show that
\begin{equation}\label{eq: nilpotentprod1}
  (x_1 + y_1)(n_1) ... (x_{2t-1} + y_{t-1})(n_{2t-1}) (x_{2t} + y_{2t}) = 0
\end{equation}
for all integers $n_1, ..., n_{t-1}$. Expanding (\ref {eq: 
nilpotentprod1}) yields a sum
of terms, each of which is a product
\begin{equation}\label{eq: nilpotentprod2}
  a_1(n_1)  ...  a_{2t-1}(n_{2t-1}) a_{2t}
\end{equation}
with each $a_i$ equal to either $x_i$ or $y_i$. So it suffices to 
show that each expression
(\ref {eq: nilpotentprod2}) is equal to zero. Now either
there are (at least) $t$ indices $i$ for which $a_i = x_i$, or $t$
indices for which $a_i = y_i$.
We will assume that the first case holds; if it does not, then the 
same argument applies with
$x_i$ replaced by $y_i$.

Using the relation
\begin{eqnarray*}
a_i(n) a_{i+1}(m) =
  a_{i+1}(m)a_i(n) + \sum_{j \geq 0} (a_i(j)a_{i+1})(m+n-j)
\end{eqnarray*}
together with the containment $a_i(j)a_{i+1} \in I^1$ whenever $a_i \in I^1$,
we see that (\ref{eq: nilpotentprod2}) can be reexpressed as a sum of terms
of the form
\begin{equation}\label{eq: nilpotentprod3}
Ab_1(m_1) ... b_k(m_k) a_{2t}
\end{equation}
where each $b_i$ is some $x_j; k \geq t-1$ and either $k \geq t$ or 
$a_{2t} = x_{2t}$;
$A$ is a product of other Fourier modes (of no interest to us). 
Because $(I^1)^t = 0$,
we have $b_1(m_1) ... b_k(m_k)a_{2t} = 0$, so that the expression 
(\ref {eq: nilpotentprod3})
vanishes. Therefore so too does (\ref{eq: nilpotentprod2}), and part 
(a) is proved.

Next, it is clear that if $V$ is semisimple then $N(V) = 0$, and that if
$ I \subseteq V$ is an ideal then $I +N(V)/I \subseteq N(V/I)$. Part (b)
is an immediate consequence of these observations.

As for (c),  after
Theorem 3.7 it is enough to show that the ideal $I$ of $V$
generated by $J' = J(Z(V))$ is nilpotent.  Because $J'$
is a nilpotent ideal in $Z(V)$, there is an integer
$t$ such that $(J')^t = 0$. We show that $I^t = 0$.
Elements of $I^t$ are sums of states of the form
\begin{equation}\label{eq: nilpotentprod4}
(a_1(-1)b_1)(n_1) ... (a_{t-1}(-1)b_{t-1})(n_{t-1}) a_t(-1)b_t
\end{equation}
with $a_1, ..., a_t \in J^{'}, b_1, ..., b_t \in V$ and $n_1, ..., 
n_{t-1} \in Z$
(cf. Lemma 2.1). Indeed, Lemma 2.1 together with associativity
also shows that $(a(-1)b)(n) = a(-1)b(n)$ for $a \in J^{'}, b \in V$ 
and $n \in Z$.
Then expression (\ref {eq: nilpotentprod4}) is equal to
\begin{eqnarray*}
a_1(-1)b_1(n_1) ... a_{t-1}(-1)b_{t-1}(n_{t-1})a_t(-1)b_t  &=& \\
  b_1(n_1) ... b_{t-1}(n_{t-1})b_t(-1) a_1(-1) ... a_t(-1)\mathbf{1} &=& 0,
\end{eqnarray*}
where we used Lemma 2.4(a) for the first equality. This completes the
proof of part (c). \qed

Using similar arguments to the above one can also establish the following
(details left to the reader): a vertex operator algebra $V$ has
a unique maximal solvable ideal $B(V)$; $B(V/B(V)) = 0$;
  $B(V) = 0$ if, and only if, $N(V) = 0$;
  $B(V) \subseteq J(V)$.  We have the following containments:
\begin{eqnarray*}
J(Z(V)) \subseteq N(V) \subseteq B(V) \subseteq J(V) \subseteq V.
\end{eqnarray*}

\section{ Proof of Theorems 1 and 2}
\setcounter{equation}{0}

   We define
  a \emph{local} or \emph{semilocal} vertex operator algebra as in the 
Introduction.
  Note that a local vertex operator
algebra is indecomposable. From the decomposition
into blocks, we see that for $V$ to be semilocal,
it is sufficient that every block is a local vertex operator
algebra. One is thus led to ask the following question.

\begin{eqnarray}\label{eq: question}
  &\mbox{Suppose that} \  V \ \mbox{is an indecomposable vertex} \\
& \mbox{operator algebra. Is it true that $V$ is local?} \nonumber
\end{eqnarray}
  An affirmative answer would reduce questions about general vertex operator
algebras to questions about local vertex operator algebras.
  Note that by Proposition 2.6, (\ref{eq: question}) is essentially 
the same as the question of
lifting idempotents of vertex operator algebras (cf. the discussion 
prior to (\ref {eq: 1decomp})),
in analogy with the familiar result from the theory of associative algebras.
\begin{eqnarray}\label{eq: idemplift}
  & \mbox{Suppose that $e \in V$ is an idempotent in  $V/J(V)$.}  \\
& \mbox{Can $e$  be lifted to an idempotent in  $V$ ?} \nonumber
   \end{eqnarray}

We now turn our attention to the proof of
Theorems 1 and 2, which will show that the answers to (\ref{eq: question})
and (\ref {eq: idemplift}) are affirmative in case
the appropriate conditions are satisfied. Our  approach
to the first part of Theorem 1 involves
the strucure of the nonassociative algebra on $V_0$
induced by the $-1$ product, which contains $Z(V)$
as a subalgebra. With (\ref {eq: idemplift}) in mind, it comes down to
showing that all idempotents in $V_0$ are contained in $Z(V)$.
We begin by developing the necessary background
concerning \emph{power associative} algebras. Albert developed this
theory in several papers (cf. [A1], [A2], [S]). It fits well into the
vertex operator algebra formalism.

A finite-dimensional nonassociative (unital) algebra $A$
is called power associative if each element $a \in A$ has the
property that all powers of $a$ associate. Thus, the subalgebra of $A$
generated by $a$ is associative and commutative. By [A1], it is sufficient
that the identities
\begin{equation}\label{eq; assocaite}
aa^2 = a^2a, \ \ \ a(aa^2) = a^2a^2
\end{equation}
hold. If $A$ is power associative, one may unambigously refer
to \emph{nilpotent} elements ($a^n=0$ for some $n \geq 1$). An
ideal $I \subseteq A$ is \emph{nil} in case every
element of $I$ is nilpotent. Then $A$ has a unique maximal
nil ideal $N$ (the nil radical), and the quotient $A/N$ has trivial
nil radical. We call $A$ a \emph{local} power associative algebra in case
its nil radical is the unique maximal ideal of $A$ and has
codimension $1$.

\begin{prop}\label{p4.1} Suppose that $A$ is a commutative, power associative
unital algebra with the property that $1$ is the only nonzero 
idempotent in $A$.
Then $A$ is a local algebra.
\end{prop}

\pf  If $A$ is also assumed to be simple, Albert proved ([A2, 
Theorem 9])
that $A = \mathbf{C}1$. So in this case the nil radical is trivial and
we are done. So we may assume that there is a nonzero proper ideal 
$I\subseteq A$.

Choose $a \in I$. The subalgebra $\langle a \rangle$ generated by $a$ 
is associative,
and hence is either nilpotent or contains a nonzero idempotent. In 
the latter case, the
hypothesis of the Theorem shows that $1 \in \langle a \rangle 
\subseteq I$. But then
$I = A$, a contradiction. Therefore $a$ is nilpotent. This shows that $I$
is a nil ideal, and hence contained in the nil radical $N$ of $A$.

We assert that the hypotheses of the Theorem hold in the
quotient algebra $A/I$. Once this is established, the Theorem follows
easily by induction on dim$A$. Now $A/I$ is certainly commutative
and power associative, so to prove the assertion it
suffices to show that
$1+I$ is the only nonzero idempotent in $A/I$. Let $e+I$ be a
nonzero idempotent
in $A/I$. It is apparent that the subalgebra $\langle e \rangle$
of $A$ generated by $e$ cannot be nilpotent. As before, this leads to
$1 \in \langle e \rangle$, and even $\langle e \rangle = \mathbf{C}1 
\oplus J(\langle e \rangle)$.
So $e = \alpha1 + x$ for some scalar $\alpha \in \mathbf{C}$ and
where $x \in  J(\langle e \rangle)$ is nilpotent. Then
$e^2 - e = (\alpha^2 - \alpha)1 + y$ with $y = x^2 + (2\alpha-1)x$ 
nilpotent. But $e^2-e \in I$ and
hence is itself nilpotent. It follows that $\alpha^2 = \alpha = 0$ 
\mbox{or} \ $1$.
As e is not nilpotent then the case $\alpha = 0$ is impossible, so 
$\alpha = 1$ and
   \begin{equation}\label{eq: twosum}
  e = 1+x; \ \ y = x^2+x \in I.
\end{equation}

Let $n \geq 1$ be the smallest positive integer satisfying $x^n \in I$.
Because $x$ is nilpotent, such an $n$ certainly exists. If $n \geq 2$ then
\begin{eqnarray*}
I = x^n+I = x^{n-2}x^2+I = -x^{n-2}x +I,
\end{eqnarray*}
where we
used the second equality in  (\ref{eq: twosum}). Then $x^{n-1} \in 
I$, contradicting
the definition of $n$. We conclude that $n=1$. Then $x \in I$
and $e+I = 1+I$, as we see from the first equality in (\ref{eq: twosum}).
This completes the proof of the Proposition. \qed

For the remainder of this Section we fix $V$ to be a
vertex operator algebra satisfying the truncation condition
\begin{equation}\label{eq: weighttruncate}
V_n = 0 \ \mbox{for} \ n \leq -2.
\end{equation}

\begin{prop}\label{p4.2} With respect to the $-1$ operation,
$V_0$ is a power associative algebra. Moreover if $a \in V_0$
then
  \begin{equation}\label{eq: vertcommut}
[a(m), a(n)]= 0 \ \mbox{for all integers} \ m,n.
  \end{equation}
\end{prop}
\pf  Pick $a \in V_0$. In addition to
(\ref {eq: vertcommut}) we have to establish the two identities
(\ref{eq; assocaite}), which amount to the following:
\begin{eqnarray}\label{eq: -powerassoc}
a(-1)^2a &=& (a(-1)a)(-1)a; \nonumber \\
a(-1)^3a &=& (a(-1)a)(-1)a(-1)a.
\end{eqnarray}
Now we have $a(n)a \in V_{-n-1}$, and in particular
\begin{equation}\label{eq: anavanish}
a(n)a =0 \ \mbox{for} \ n \geq 1
\end{equation}
  thanks to (\ref {eq: weighttruncate}). By
skew-symmetry we also have
\begin{eqnarray*}
a(0)a = -a(0)a - \sum_{n \geq 1} (-1)^n \frac{L(-1)^n}{n!}a(n)a.
\end{eqnarray*}
Then by (\ref {eq: anavanish}) we also get $a(0)a = 0$.
So in fact $a(n)a = 0$ for all $n \geq 0$, and (\ref {eq: vertcommut})
is a consequence of this together with the commutivity axiom.
  With (\ref {eq: vertcommut}) in hand, the proof of (\ref {eq: -powerassoc})
is a straightforward application of the associativity axiom. We omit 
the details.
\qed

  The next result contains a key calculation.

\begin{prop}\label{p4.3} $Z(V)$
contains every idempotent of $V_0$.
\end{prop}
\pf  We begin with the observation that the idempotents and
  the involutorial units of $V_0$ (i.e., those $u \in V_0$ satisfying
$u(-1)u = \mathbf{1}$) span the \emph{same} subspace.
This is a simple fact about any unital nonassociative algebra,
  and we omit the proof. In view of this, in order to prove the
Proposition it suffices to show that $Z(V)$ contains
every involutorial unit $u \in V_0$. Fix such an element $u$.
We must show that $L(-1)u = 0$. Since $L(-1)\mathbf{1} = 0$ then
\begin{eqnarray*}
0 &=& u(-1)L(-1)u(-1)u \\
     &=& u(-1)( [L(-1), u(-1)]u + u(-1)L(-1)u(-1)\bf{1} )\\
     &=&u(-1) (u(-2)u + u(-1)[L(-1), u(-1)] \ \bf{1}) \\
      &=& u(-1) (u(-2)u + u(-1)u(-2)\ \bf{1}) \\
    &=& 2u(-1)u(-2)u\\
    &=& 2u(-2)\bf{1} \\
   &=& 2(L(-1)u)(-1)\bf{1} \\
   &=&  2L(-1)u,
\end{eqnarray*}
where we have made use of  (\ref {eq: vertcommut}).    \qed

We turn to the proof of Theorem 2. That (a) $\Rightarrow$ (b) is clear, while
(b) $\Leftrightarrow$ (c) is nothing but
  Proposition 2.6.
  If (c) holds
then by the previous Proposition, $\mathbf{1}$
is the unique nonzero idempotent in $V_0$.
Let $V_0^+$ be the algebra obtained from
$V_0$ by redefining the product to be
$a*b = 1/2(a(-1)b + b(-1)a)$. Then
$V_0^+$ is a unital commutative algebra, and
it is power associative because $V_0$ is.
Furthermore, the idempotents of $V_0$
and $V_0^+$ are the same, whence $\mathbf{1}$ is also
the unique nonzero idempotent of $V_0^+$. By Proposition 4.1,
$V_0^+$ is a local algebra, and so (d) holds.

To complete the proof of Theorem 2, we must show that (d) $\Rightarrow$ (a).
Suppose then, that $V_0^+$ is a local power associative algebra. Any
proper ideal $I \subseteq V_0$ is an ideal in $V_0^+$ and therefore
contained in the nil radical of $V_0^+$. In particular,
$V_0^+/I$ is again a local power associative algebra with the identity
element the only nonzero idempotent. It follows from this that $V_0/I$
also has a unique nonzero idempotent.

Set $J = J(V)$. Being semisimple,
$V/J$ is the sum of $n \geq 1$ simple vertex operator algebras.
   Choose states $e_i \in V_0$
such that $e_i + J = \mathbf{1}_i  \ (1 \leq i \leq n)$ are the vacuum states
of the simple components of $V/J$. Now $J_0$ is an ideal of $V_0$
and evidently the $e_i$ map onto \emph{distinct} idempotents of $V_0/J_0$.
After the last paragraph, we conclude that in fact $n=1$, that is $V$ 
is a local
vertex operator algebra. This completes the proof of Theorem 2. \qed

Theorem 1(i) is now an easy consequence of what we have already established.
Namely, since $V$ satisfies the condition (\ref {eq: weighttruncate}) it is
clear that the blocks
  $V^i$ of $V$ have the same property. Since the blocks are indecomposable,
  Theorem 2 informs us that they are local.
This completes the proof of Theorem 1(i). \qed

We turn to the proof of the second part of Theorem 1. So we are assuming that
the vertex operator algebra $V$ is such that the Jacobson radical 
coincides with the solvable radical.
We proceed, using induction on the dimension of $V_0 + V_1$, to show that
$V$ is semilocal. The decomposition into blocks shows that we may assume
without loss that $V$ is indecomposable. Thus, $Z(V)$ is a local
algebra by Proposition 2.6, and we must show that $V$ is a local 
vertex operator algebra.
  If $J = J(V)$ is zero then the result
is clear, so we may take $J \neq 0$. We may also assume that $T(V) = 0$
by Lemma 3.1, in which case $J$ contains a minimal ideal $N$ by
Lemma 3.2. Since $J$ is solvable then $N$ is necessarily
nilpotent, indeed $N^2 = 0$.

Set $A =$ Ann$_V(N)$. Then $A$ is an ideal of $V$ which contains $N$, and
$V/A$ is indecomposable by Lemma 2.8. On the other hand, by induction we
know that $V/N$ is semilocal. Let
\begin{eqnarray*}
V/N = U^1/N \oplus ... \oplus U^n/N
\end{eqnarray*}
be the decomposition of $V/N$ into blocks. Thus, each $U^i/N$ is a local
vertex operator algebra. Pick states $a_1, ... , a_n \in V_0$ such that
$a_i + N$ is the vacuum state of $U^i/N$. We claim that all but one 
of the states
$a_1, ... , a_n$ lie in $A$. If $n=1$ this is clear, so assume that 
$n \geq 2$ and
that neither $a = a_1$ nor $b = a_2$ lie in $A$. In this case both $a$ and
$b$ map onto the vacuum of $V/A$ (because it is indecomposable)
and therefore $a - b \in A$.  So $a - b \in A \cap (U^1+U^2)$.
Since each $U^i$ is local it follows that $U^1+U^2 = J(U^1+U^2) + A 
\cap (U^1+U^2)$,
so that $a, b \in U^1+U^2 = A \cap (U^1+U^2) \subseteq A$, a contradiction.
So indeed we may assume $a_1, ... , a_{n-1}$ are contained in $A$. Now if $n=1$
then there is nothing to prove, so we may assume without loss
  that $n \geq 2$ and $a = a_1 \in A$.
We will show that this leads to a contradiction.

Since $a+N$ is the vacuum state of $U^1/N$ then $a(-1)a = a + n$
for some $n \in N_0$. Because $a \in A$ then $a(k)N = 0$ for all integers $k$,
and from this together with $N^2 = 0$
  it follows easily that $(a+n)(-1)(a+n) = a+n$. So without loss
we may assume that $a(-1)a = a$. We will show that $a \in Z(V)$. In this case
$a$ is a nonzero idempotent in $Z(V)$ and hence $a = \mathbf{1}$ because $Z(V)$
is a local algebra. Then $\mathbf{1} \in A$ and therefore annihilates 
$N$, contradiction.

The strategy used to establish $a \in Z(V)$ is similar to that used in
  the proof of Proposition 4.3. First note that since $a+N$ is the vacuum of
$U^1/N$ then $a(k)V \subseteq N$ for $k \neq -1$. In particular, 
because $a \in A$ we see that
$a(-1)$ annihilates $a(k)a$ and $(a(k)a)(l)$ annihilates $a$  for $k 
\neq -1$ and all $l$.
 From this we conclude that
\begin{eqnarray}\label{eq: kill1}
a(-2)a &=& a(-2)a(-1)a \nonumber \\
           &=& a(-1)a(-2)a + [a(-2), a(-1)]a \nonumber \\
&=&a(-1)a(-2)a +  \sum_{i \geq 0} \left( \begin{array}{c}
                                               -2 \\  i
                                         \end{array}  \right) 
(a(i)a)(-3-i)a \nonumber \\
  &=& 0.
  \end{eqnarray}
Similarly we get
\begin{equation}\label{eq: kill2}
a(-1)a(-2)\mathbf{1} = 0.
\end{equation}
  Now apply the operator $L(-1)$ to both sides of the equality
$a(-1)a = a$ to see that
\begin{eqnarray*}
L(-1)a  &=& L(-1)a(-1)a \\
  &=&  [L(-1), a(-1)]a + a(-1)L(-1)a \\
&=& a(-2)a + a(-1)a(-2)a\\
&=& 0, \\
  \end{eqnarray*}
  where we used (\ref{eq: kill1}) and (\ref{eq: kill2}).
This shows that indeed $a \in Z(V)$, and the proof of the second part 
of Theorem 1 is complete.
\qed

\section{Examples}
\setcounter{equation}{0}

We discuss some examples of Jacobson radicals and local vertex 
operator algebras.
Simple vertex operator algebras are obvious examples of local vertex
operator algebras.

{\sc Example 1:} Semidirect products.\\
Suppose that $M$ is a $Z$-graded module over a vertex operator algebra $V$.
In [L], Li describes how to construct the semidirect product
\begin{equation}\label{eq: dirsum}
W = V \oplus M,
  \end{equation}
where in particular $V$ is a vertex operator subalgebra of $W$,
$M$ is an ideal,  and
  $Y_W(u, z)v = 0$ for all $u, v \in M$. So $M$ is a nilpotent ideal
of $W$, and it is evident that if $V$
is semisimple then $M = J(W) = N(W)$. In particular, if $V$ is simple then
$W$ is a local vertex operator algebra.

{\sc Example 2:}
Let $V_L$ be the lattice vertex operator algebra
associated with the root lattice of type $A_1$ spanned by a single 
positive root $\alpha$, say.
  Let $M(1) \subseteq V_L$
be the Heisenberg vertex operator 
subalgebra. 
  We have the decomposition into simple $M(1)$-modules
\begin{eqnarray*}
V_L = \bigoplus_{n \in Z}M(1) \otimes \mathbf{C}e^{n \alpha}.
\end{eqnarray*}
Set
\begin{eqnarray*}
U = \bigoplus_{n \geq 0}M(1) \otimes \mathbf{C}e^{n \alpha} .
  \end{eqnarray*}
Then $V_L$ is simple and $U$ is a subvertex operator algebra of $V_L$
which is also a local vertex operator algebra with Jacobson radical
\begin{eqnarray*}
J(U) =  \bigoplus_{n > 0}M(1)\otimes \mathbf{C}e^{n \alpha} .
\end{eqnarray*}
This shows that the Jacobson radical is not \emph{functorial}. That is,
if $f: V_1 \rightarrow V_2$ is a morphism of vertex operator algebras, then
  $f$ does \emph{not} necessarily induce a morphism $f: J(V_1) 
\rightarrow J(V_2)$.

  {\sc Example 3: Virasoro vertex operator algebras}.\\
  Consider
the Virasoro vertex operator algebra $V_c$ of central charge $c$
defined as the quotient of the Verma module over the Virasoro
algebra $Vir$
of highest weight zero and central charge $c$ by the submodule generated
by $L(-1)\bf{1}$. (See [FZ], [W] for more details.) If we choose $c$ to be in
the \emph{discrete series} (loc. cit.) then $V_c$
has a submodule $M \neq 0$ such that
\begin{eqnarray*}
V_c/M = L_c
\end{eqnarray*}
  is the simple Virasoro
vertex operator algebra $L_c$. Clearly $V_c$ is a
local vertex operator algebra and $J(V_c) =M$. Note that
a decomposition $V_c = M \oplus L_c$ with $L_c$
a vertex operator subalgebra
is impossible,  because $V_c$ is generated by $\bf{1}$ as a module over $Vir$.

\end{document}